\documentclass[reqno,12pt]{amsart}
\usepackage{amsmath, latexsym, amsfonts, amssymb, amsthm, amscd}

\numberwithin{equation}{section}

\newtheorem{theorem}{Theorem}[section]
\newtheorem{lemma}[theorem]{Lemma}
\newtheorem{proposition}[theorem]{Proposition}
\newtheorem{corollary}[theorem]{Corollary}
\newtheorem{rem}[theorem]{Remark}

\renewcommand{\tilde}{\widetilde}          
\DeclareMathSymbol{\leqslant}{\mathalpha}{AMSa}{"36} 
\DeclareMathSymbol{\geqslant}{\mathalpha}{AMSa}{"3E} 
\DeclareMathSymbol{\eset}{\mathalpha}{AMSb}{"3F}     
\renewcommand{\leq}{\;\leqslant\;}                   
\renewcommand{\geq}{\;\geqslant\;}                   

\newcommand{\suptwo}[2]{\sup_{\substack{#1 \\ #2}}} 
\newcommand{\inftwo}[2]{\inf_{\substack{#1 \\ #2}}} 
\newcommand{\sumtwo}[2]{\sum_{\substack{#1 \\ #2}}} 

\newcommand{\R}{\mathbb{R}}
\newcommand{\Z}{\mathbb{Z}}
\newcommand{\N}{\mathbb{N}}

\title[A local limit theorem for directed polymers]{A local limit
  theorem  for directed polymers in random media: the continuous and the discrete case}

\author{Vincent Vargas}
\thanks{Partially
   supported by CNRS (UMR 7599
``Probabilit{\'e}s et Mod{\`e}les
Al{\'e}atoires'')}
        
\begin{document}

\maketitle
\begin{center}
{\footnotesize \noindent
 Universit{\'e} Paris 7,\\
Math{\'e}matiques, case 7012,\\ 2, place Jussieu, 75251 Paris, France}

{\footnotesize \noindent e-mail: \texttt{vargas@math.jussieu.fr}}
\end{center}

\begin{abstract}
In this article, we consider two models of directed polymers in random environment: a
discrete model in a general random environment and a continuous model. We consider these models in
dimension greater or equal to 3 and we suppose that the normalized
partition function is bounded in $L^{2}$ (the "high" temperature case).
Under these
assumptions, Sinai proved in \cite{cf:Sin} a local limit theorem for
the discrete model, using a perturbation expansion. In this article,
we give a new method for proving Sinai's local limit theorem. This new
method can be transposed to the continuous setting in which
we prove a similar local limit theorem.

\bigskip

\noindent\textsc{Resum{\'e}:} Dans  cet article, on consid{\`e}re deux mod{\`e}les
de polym{\`e}res dirig{\'e}s en environnement al{\'e}atoire: un mod{\`e}le
discret en environnement al{\'e}atoire g{\'e}n{\'e}ral et un mod{\`e}le continu. On
consid{\`e}re ces mod{\`e}les en dimension sup{\'e}rieur ou {\'e}gale {\`a} 3 et on
suppose que la fonction de partition renormalis{\'e}e est born{\'e}e dans
$L^{2}$ (cela correspond au cas de "haute" temp{\'e}rature). Sous ces
hypoth{\`e}ses, Sinai a montr{\'e} dans \cite{cf:Sin} un th{\'e}or{\`e}me limite
locale pour le mod{\`e}le discret en utilisant un d{\'e}veloppement en
perturbation. Dans cet article, on donne une nouvelle m{\'e}thode pour
d{\'e}montrer le th{\'e}or{\`e}me limite locale ci-dessus. Cette nouvelle m{\'e}thode
peut {\^e}tre transpos{\'e}e au cas continu dans lequel on montre un th{\'e}or{\`e}me
limite locale similaire.

\bigskip


\bigskip

 
\end{abstract}

\section{Introduction}
Directed polymers in random environment is a model of
statistical mechanics in which stochastic processes interact with a
random environment, depending on both time and space: one studies the
path of the stochastic process under a random Gibbs measure depending
on the temperature (as the temperature increases, the influence of the
random environment decreases).

In this article, we will consider two polymer models: a simple random walk
model of directed polymers and its continuous analogue, a Brownian
model of directed polymers. The discrete model first
appeared in the physics litterature (\cite{cf:HuHe85}) to modelize the
phase boundary of Ising model subject to random impurities and its
first mathematical study was undertaken by Imbrie, Spencer in 1988 
(\cite{cf:ImSp88}) and Bolthausen in 1989 (\cite{cf:Bol89}). The
continuous model we study here was first introduced and studied by
Comets and Yoshida in 2004 (\cite{cf:CY03}). These models are related
to many models of statistical physics. We refer to the
survey paper \cite{cf:KrSp91} by Krug and Spohn  for an account on these
models and there relations.    

In the sequel, we will suppose that the dimension of the underlying
stochastic process is greater or equal
to 3 and that the normalized partition function is bounded in $L^{2}$
(see subchapters 1.1-1.2. for the definition of the normalized
partition function). Under these assumptions, the polymer is
diffusive in the sense that a central limit theorem holds:  by scaling
by the square root of time, the simple random walk (Brownian motion in
the continuous model) converges in law under the random Gibbs measure
to a gaussian measure (see \cite{cf:ImSp88}, \cite{cf:Bol89}, \cite{cf:SoZh96},
\cite{cf:CY04}). One can sometimes go a step
further than convergence in law by
giving an equivalent of the density: this is called a local limit
theorem. In \cite{cf:Sin}, Sinai obtained a local limit theorem by
using a perturbation expansion. Unfortunately, it is not clear how to
adapt the strategy
to the continuous setting. The object of this work is to give a
new method for proving Sinai's theorem; this method is sufficiently
general to be easily adapted to prove a similar local limit theorem
in the continuous setting. Our approach is simple and relies only on
$L^2$ computations and on properties of the simple random walk bridges
(Brownian bridges in the continuous case).

Finally, we remind that some results have been achieved in the case
of dimension less or equal to 2 or when the temperature is low. In
these cases, the polymer is non-diffusive (see remark \ref{rem:Com} below) and many conjectures remain
open. For an account on these cases, we refer to \cite{cf:CaHu02} in a
gaussian environment and to \cite{cf:CSY04} in a general
environment.

The article is organized as follows: each chapter is divided into two
subchapters, one of them being devoted to the discrete model and the other one
being devoted to the continuous model. First, we introduce the two
models. In the second chapter, we will remind the known results 
at high temperature when the dimension of the underlying process is
greater or equal to $3$; we will also formulate an analogue to Sinai's
local limit theorem for the Brownian directed polymer. In the third chapter, we
will prove the local limit theorem for both models.

\subsection{The simple random walk model of directed polymers}

\begin{itemize}
\item
Let $((\omega_{n})_{n \in \N},(P^x)_{x \in \Z^d})$ denote the simple
random walk on the $d$-dimensional integer lattice $\Z^d$, defined on
a probability space $(\Omega, \mathcal{F})$; more precisely, for $x$
in $\Z^d$, under the measure $P^x$, $(\omega_{n}-\omega_{n-1})_{n \geq
  1}$ are independant and 
\begin{equation*}
P^x(\omega_{0}=x)=1, \; \; P^x(\omega_{n}-\omega_{n-1}=\pm
\delta_{j})=\frac{1}{2d},\; \; j=1, \ldots, d, 
\end{equation*}
where $(\delta_{j})_{1 \leq j \leq d}$ is the j-th vector of the
canonical basis of $\Z^d$. In the sequel, $P$ will denote $P^0$.
For $x$ in $\Z^d$, let $q^{(n)}(x)$ be the
probability for the random walk starting in $0$ to be in $x$ at time
$n$:
$$q^{(n)}(x)\overset{def.}{=}P(\omega_{n}=x).$$

\item
 The random environment
on each lattice site is a sequence $\eta=(\eta(n,x))_{(n,x) \in \N
  \times \Z^d}$ of real valued, non-constant and
i.i.d. random variables defined on a probability space $(H,
\mathcal{G}, Q)$ such that
$$\forall \beta \in \R \; \; \lambda(\beta)\overset{def.}{=} \ln Q(e^{\beta  \eta(n,x)})< \infty.$$
\item
For any $n>0$, we define the (Q-random) polymer measure $\mu_{n}^x$ on
the path space $(\Omega, \mathcal{F},P^x)$ by:
$$\mu_{n}^x(d\omega)=\frac{1}{Z_{n}^x}\exp(\beta H_{n}(\omega)-n\lambda(\beta))P^x(d\omega)$$
where $\beta \in \R$ is the inverse temperature, 
$$
 H_{n}(\omega)\overset{def.}{=}\sum_{j=1}^{n}\eta(j,\omega_{j}) \qquad
$$
and
$$
Z_{n}^{x}=P^x(\beta H_{n}(\omega)-n\lambda(\beta)) \qquad
$$
is the normalized partition function ($Q(Z_{n}^{x})=1$).
\end{itemize}
Let $(\mathcal{G}_{n})_{n \geq 0}$ be the filtration defined by

$$\mathcal{G}_{n}=\sigma \lbrace \eta(j,x); j \leq n, x \in \Z^d
\rbrace.$$
For any fixed path $\omega$,
$((\sum_{j=1}^{n}\beta\eta(j,\omega_{j}))-n\lambda(\beta))_{n \geq 1}$ is
a random walk with independent increments thus it is not hard to see
that $(Z_{n}^{x},\mathcal{G}_{n})_{n \geq 0}$ is a positive
martingale. Therefore, it converges $Q$-a.s. to a limit $Z_{\infty}^x$.
Since the event $(Z_{\infty}^x=0)$ is measurable with respect to the
tail $\sigma$-field
$$\underset{n \geq1}{\bigcap}\sigma \lbrace \eta(j,x); j \geq n, x \in \Z^d
\rbrace,$$
by Kolmogorov's $0-1$ law, there are only two possible situations 
$$Q(Z_{\infty}^x=0)=1 \; \; \; or \; \; \; Q(Z_{\infty}^x=0)=0.$$

In the former case, we say that strong disorder holds and in the
latter case we say that weak disorder holds.

\subsection{The Brownian motion model of directed polymers}

\begin{itemize}
\item
Let $((\omega_{t})_{t \in \R_{+}},(P^x)_{x \in \R^d})$ denote a
$d$-dimensional standard brownian motion, defined on a probability
space $(\Omega, \mathcal{F})$. In the sequel, $P$ will denote $P^0$. 
For $t > 0$ and $x,y$ in $\R^d$, 
let $p(t,x,y)$ be the transition density of the Brownian motion:
$$p(t,x,y)\overset{def.}{=}\frac{1}{(2\pi t)^{\frac{d}{2}}}e^{-\frac{\mid y-x
    \mid^2}{2t}}.$$

\item
The random environment $\eta$ is a Poisson
random measure on $\R_{+} \times \R^d$ with unit intensity, defined on a
probability space $(M,\mathcal{G}, Q)$. We recall that $\eta$ is an
integer valued random measure characterized by the following property:
If $A_{1}, \ldots, A_{n} \in \mathcal{B}(\R_{+} \times \R^d)$ are
disjoint and bounded Borel sets, then
$$Q(\bigcap_{j=1}^{n}(\eta(A_{j})=k_{j}))=\prod_{j=1}^{n}e^{-\mid A_{j}\mid}\frac{\mid A_{j}\mid^{k_{j}}}{k_{j}!}.$$
where $k_{1}, \ldots k_{n} \in \N$ and $\mid. \mid$ denotes the
Lebesgue measure in $\R^{d+1}$. We define $V_{t}$ to be the unit
volume tube around the graph $\lbrace (s,\omega_{s})\rbrace_{0
  < s \leq t}$ of the Brownian path:
$$V_{t}=V_{t}(\omega)=\lbrace (s,x); s \in ]0,t], x \in
U(\omega_{s})\rbrace$$
 where $U(x)$ is the closed ball in $\R^{d}$ with unit volume and
 centered at $x \in \R^{d}$.

\item 
For any $t>0$, we define the ($Q$-random) polymer measure
$\mu_{t}^{x}$ on the path space $(\Omega, \mathcal{F},P^x)$ by:
$$\mu_{t}^{x}(d\omega)=\frac{\exp(\beta
    \eta(V_{t})-\lambda(\beta)t)}{Z_{t}^{x}}P^{x}(d\omega),$$
where $\beta \in \R$ is the inverse temperature and 
$$Z_{t}^{x}=P^{x}(\exp(\beta
    \eta(V_{t})-\lambda(\beta)t))$$
is the normalized partition function $(Q(Z_{t}^{x})=1)$. In this
    setting, the random environment is a Poisson point process so we
    get the explicit value:
$$\lambda(\beta)=e^{\beta}-1 \in ]-1, \infty[.$$
\end{itemize}
It is natural to introduce the filtration
$(\mathcal{G}_{t})_{t > 0}$ defined by :

$$\mathcal{G}_{t}=\sigma \lbrace \eta(A); A \in \mathcal{B}(]0,t]
\times \R^d) \rbrace.$$

As in the discrete setting, it is not hard to show that
$(Z_{t}^{x},\mathcal{G}_{t})_{t > 0}$ is a positive martingale wich
converges $Q$-a.s. to a non negative random variable $Z_{\infty}^{x}$
that has the following property:
 
$$Q(Z_{\infty}^x=0)=1 \; \; \; or \; \; \; Q(Z_{\infty}^x=0)=0.$$

In the former case, we say that strong disorder holds and in the
latter case we say that weak disorder holds.

\section{Study of the directed polymers when the normalized partition
  function is bounded in $L^2(Q)$} 

From now on, in the rest of this paper, we will only consider the case
$d \geq 3$ and we will suppose that that the normalized partition
function is bounded in $L^2(Q)$. In that case, the latter converges
$Q$-a.s. and in $L^2(Q)$ to the random variable $Z_{\infty}^{x}$. The
$L^2$-convergence implies that $Q(Z_{\infty}^{x})=1$ and therefore
weak disorder holds. Under these assumptions, the behavior of the
typical path under the
polymer measure is diffusive (see \cite{cf:CSY04} for the discrete case and
 \cite{cf:CY04} for the continuous case).

\subsection{The Simple random walk model}

In order to get a nice probabilistic interpretation, we work on the
product space $(\Omega^2,\mathcal{F}^{\otimes 2},(P^{x} \otimes
P^{y})_{x,y \in \Z^d})$ and thus consider another simple random walk
$(\tilde{\omega}_{n})_{n \in \N}$ independant of the first one
$(\omega_{n})_{n \in \N}$ under the same environment.

Let
$\lambda_{2}(\beta)\overset{def.}{=}\lambda(2\beta)-2\lambda(\beta)$
and $N_{k,n}=N_{k,n}(\omega,\tilde{\omega})$ be the number of ordered intersections
of $\omega$ and $\tilde{\omega}$ between $k$ and $n$:
$$N_{k,n}\overset{def.}{=}\sum_{j=k}^{n}1_{\omega_{j}=\tilde{\omega}_{j}}.$$
With these notations, the following proposition is straightforward
(e.g., \cite{cf:CSY04}): 
\begin{proposition}
We have the following identity:
\begin{align*}
Q((Z_{n}^{x})^2)&=P^{x} \otimes
P^{x}(e^{\lambda_{2}(\beta)N_{1,n}}) \\
&=P \otimes
P(e^{\lambda_{2}(\beta)N_{1,n}}).
\end{align*}
In particular,
$$\sup_{n \geq 0}Q((Z_{n}^{x})^2)=P \otimes
P(e^{\lambda_{2}(\beta)N_{1,\infty}}).$$
We have the following equivalence  
$$P \otimes
P(e^{\lambda_{2}(\beta)N_{1,\infty}}) <
\infty \; \; \Longleftrightarrow \; \; \lambda_{2}(\beta)<\ln(\frac{1}{\pi_{d}})$$
where $\pi_{d}\overset{def.}{=}P(\exists n \geq 1, \; \omega_{n}=0) < 1$.
Thus, We have the following equivalence:

$$\sup_{n \geq 0}Q((Z_{n}^{x})^2) < \infty \; \; \Longleftrightarrow \; \; \lambda_{2}(\beta)<\ln(\frac{1}{\pi_{d}}).$$

\end{proposition}
 
A serie of articles \cite{cf:ImSp88}, \cite{cf:Bol89}, \cite{cf:SoZh96}
lead to the following central limit theorem:

\begin{theorem}[Central limit Theorem]\label{th:clt}
Suppose that the normalized partition function is bounded in $L^2$:
$$\lambda_{2}(\beta)<\ln(\frac{1}{\pi_{d}}).$$
Then, for all $f \in C(\R^d)$ with at most polynomial growth at
infinity,
$$\mu_{n}^x\left(f(\frac{\omega_{n}}{\sqrt{n}})\right) \underset{n \to \infty}{\longrightarrow}
\frac{1}{(2\pi)^{\frac{d}{2}}}\int_{\R^d}f(\frac{x}{\sqrt{d}})e^{-\frac{\mid
    x \mid^2}{2}}dx, \; Q-a.s.$$
\end{theorem}

A step further is to try and prove a local limit theorem: one wants to
obtain an expansion of the density $P^{x}(e_{1,n}1_{
\omega_{n}=y})$. As mentioned in the introduction, this has been
done in \cite{cf:Sin} by Sinai. In this
paper, we will give a different proof of the local limit theorem which
can be adapted to prove a continuous analogue in the Brownian
setting.   

Let us introduce a few notations that we will use in the rest of this
paper. We define for $k \leq n$
$$e_{k,n} \overset{def.}{=}\exp((\sum_{j=k}^{n}\beta
  \eta(j,\omega_{j}))-(n-k+1)\lambda(\beta))$$
and the time reversed analogue
$$\overset{\leftarrow}{e}_{k,n}\overset{def.}{=}\exp((\sum_{j=k-1}^{n-1}\beta
  \eta(n-j,\omega_{j}))-(n-k+1)\lambda(\beta)).$$
We can now recall Sinai's local limit theorem in a suitable form:

 \begin{theorem}[Sinai, 1995]\label{th:sin}
Let $d \geq 3$, $A>0$ and $\beta$ be such that $\lambda_{2}(\beta)< \ln(\frac{1}{\pi_{d}})$.
Then, if $(l_{n})_{n \geq 0}$ is a sequence of integers that tend to
infinity such that $l_{n}=o(n^a)$ with
$a<\frac{1}{2}$,   
\begin{equation} \label{eq:sin}
P^{x}(e_{1,n}\mid
\omega_{n}=y)=P^{x}(e_{1,l_{n}})P^{y}(\overset{\leftarrow}{e}_{1,l_{n}})+\delta_{n}^{x,y}
\end{equation}
with 
$$\sup_{\mid y-x \mid \leq A\sqrt{n}}Q(\mid \delta_{n}^{x,y}
\mid^2)\underset{n \to \infty}{\rightarrow} 0.$$

This leads to the following formulation that can be found in Sinai's article:
\begin{equation}\label{eq:sin2}
P^{x}(e_{1,n}\mid
\omega_{n}=y)=Z_{\infty}^{x}P^{y}(\overset{\leftarrow}{e}_{1,n})+\bar{\delta}_{n}^{x,y}
\end{equation}
with 
$$\sup_{\mid y-x \mid \leq A\sqrt{n}}Q(\mid \bar{\delta}_{n}^{x,y}
\mid)\underset{n \to \infty}{\rightarrow} 0.$$
\end{theorem}

\begin{rem} \label{rem:intui}
Intuitively, the local limit theorem asserts that, conditionnaly to
the event $(\omega_{n}=y)$, the polymer only "feels" the environment
at times  $k$ small where it stays near $x$ and at times $k$ close to
$n$ where it stays near $y$. In between , the polymer behaves like a 
conditionned simple random walk.  
\end{rem}

\begin{rem}\label{rem:sin}
 Theorem \ref{th:sin} leads to a weak form of theorem \ref{th:clt}: for all $f \in
C(\R^d)$ with compact support,
\begin{equation*}
\mu_{n}^x\left(f(\frac{\omega_{n}}{\sqrt{n}})\right) \overset{Q-Proba.}{\underset{n \to \infty}{\longrightarrow}}
\frac{1}{(2\pi)^{\frac{d}{2}}}\int_{\R^d}f(\frac{x}{\sqrt{d}})e^{-\frac{\mid
    x \mid^2}{2}}dx.
\end{equation*}
This derivation can be found in \cite{cf:Sin}.
\end{rem}
\begin{rem}\label{rem:Com}
At a heuristic level, we argue that the local limit theorem is a
natural definition for the polymer to be diffusive (more natural than
the central limit theorem itself). Roughly, the local limit theorem implies
\begin{align*} 
I_{n} &\overset{def.}{=} \sum_{x \in \Z^d} \mu_{n}(\omega_{n}=x)^2 \\ 
&\approx \sum_{x \in \Z^d} (P^{x}(\overset{\leftarrow}{e}_{1,n}))^{2}q^{(n)}(x)^2\\
&\approx Q(Z_{n}^2) \times \sum_{x \in \Z^d}q^{(n)}(x)^2\\
& \approx \frac{C}{n^{d/2}}.
\end{align*}
With other respects, recall (e.g. \cite{cf:CSY04}) that for $d=1,2$
and $\beta \not=0$ or $d \geq 3$ and $\beta$ large, 
\begin{equation*}
\exists \delta > 0, \; \; \; \varlimsup_{n \to \infty}I_{n} \geq \delta \; \; Q-a.s.
\end{equation*}
(at least if $\eta$ is unbounded in the second case). Therefore, it is
natural to call these two cases "non-diffusive" as mentionned in the introduction.

\end{rem}

\subsection{The Brownian model}
This subchapter is the continuous analogue of the previous one. We work on the
product space $(\Omega^2,\mathcal{F}^{\otimes 2},(P^{x} \otimes
P^{y})_{x,y \in \R^d})$ and thus consider another d-dimensional brownian motion
$(\tilde{\omega}_{t})_{t \in \R_{+}}$ independant of the first one
$(\omega_{t})_{t \in \R_{+}}$ under the same environment.

Let
$\lambda_{2}(\beta)\overset{def}{=}\lambda(2\beta)-2\lambda(\beta)$
where we recall that $\lambda(\beta)=e^{\beta}-1$.
Let $N_{s,t}=N_{s,t}(\omega,\tilde{\omega})$ be the volume of the overlap in
time $[s,t]$ of unit tubes around
 $\omega$ and $\tilde{\omega}$:
$$N_{s,t}\overset{def.}{=}\int_{s}^{t}\mid
U(\omega_{u})\cap U(\tilde{\omega}_{u})\mid du.$$
With these notations, we can find the following proposition in \cite{cf:CY03}: 

\begin{proposition}
We have the following identity:
\begin{align*}
Q((Z_{t}^{x})^2)&=P^{x} \otimes
P^{x}(e^{\lambda_{2}(\beta)N_{0,t}})\\
&=P \otimes
P(e^{\lambda_{2}(\beta)N_{0,t}}).
\end{align*}
In particular,
$$\sup_{t \geq 0}Q((Z_{t}^{x})^2)=P \otimes
P(e^{\lambda_{2}(\beta)N_{0,\infty}}).$$
There exists $\lambda(d) > 0$ such that:
$$ \lambda' \in ]0,\lambda(d)[ \; \; \Longleftrightarrow \; \; P \otimes
P(e^{\lambda'N_{0,\infty}}) < \infty.$$
\end{proposition}
In \cite{cf:CY04}, Comets and Yoshida prove the following central limit theorem:\begin{theorem}[Central limit theorem]
Suppose that $\beta$ is such that:
$$\lambda_{2}(\beta)<\lambda(d).$$
Then, for all $f \in C(\R^d)$ with at most polynomial growth at
infinity,
$$\mu_{t}^x\left(f(\frac{\omega_{t}}{\sqrt{t}})\right) \underset{t \to \infty}{\longrightarrow}
\frac{1}{(2\pi)^{\frac{d}{2}}}\int_{\R^d}f(\frac{x}{\sqrt{d}})e^{-\frac{\mid
    x \mid^2}{2}}dx, \; Q-a.s.$$
\end{theorem}
As in the discrete setting ,we define for $s \leq t$
$$e_{s,t} \overset{def.}{=}e^{\beta \eta(V_{s,t})-\lambda(\beta)(t-s)}$$
where $V_{s,t}$ is the unit tube around the graph
$\lbrace(u,\omega_{u})_{s < u \leq t}\rbrace$:
$$V_{s,t}=\lbrace(u,x); u \in ]s,t], x \in U(\omega_{u}) \rbrace.$$ 
We also define the time reversed analogue:
$$\overset{\leftarrow}{e}_{s,t}\overset{def}{=}e^{\beta \eta(\overset{\leftarrow}{V}_{s,t})-\lambda(\beta)(t-s)}$$
where 
$$\overset{\leftarrow}{V}_{s,t}=\lbrace(t-u,x); u \in ]s,t], x \in U(\omega_{u}) \rbrace.$$
We can now formulate a new result: the local limit
theorem for Brownian polymers.  
\begin{theorem}\label{th:sin2}
Let $d \geq 3$, $A>0$ and $\beta$ be such that $\lambda_{2}(\beta)<
\lambda(d)$. Then, if $(l_{t})_{t \geq 0}$ is a positive function that tends to
infinity such that $l_{t}=o(t^a)$ with
$a<\frac{1}{2}$,   
$$P^{x}(e_{0,t}\mid \omega_{t}=y)=P^{x}(e_{0,l_{t}})P^{y}(\overset{\leftarrow}{e}_{0,l_{t}})+\delta_{t}^{x,y}$$
with 
$$\sup_{\mid y-x \mid \leq A\sqrt{t}}Q(\mid \delta_{t}^{x,y}
\mid^2)\underset{t \to \infty}{\rightarrow} 0.$$
This leads to the following formulation in $L^1$:
$$P^{x}(e_{0,t}\mid \omega_{t}=y)=Z_{\infty}^{x}P^{y}(\overset{\leftarrow}{e}_{0,t})+\bar{\delta}_{t}^{x,y}$$
with 
$$\sup_{\mid y-x\mid \leq A\sqrt{t}}Q(\mid \bar{\delta}_{t}^{x,y}
\mid)\underset{t \to \infty}{\rightarrow} 0.$$
\end{theorem}
\begin{rem}
The remarks \ref{rem:intui} and \ref{rem:sin} apply here too.
\end{rem}

\section{Proofs}
Our proof of theorem \ref{th:sin} is based on the way bridge measures
of the simple random walk relate to the measure of the simple random
walk. This proof can be translated in the continuous
setting because Brownian bridge measures relate to the
Wiener measure in a similar way. The two main relations we use are the
absolute continuity result (\ref{eq:moi}) (relation (\ref{eq:moi'}) in the
Brownian setting) and the inequality (\ref{eq:moi2}) (relation (\ref{eq:moi2'}) in the
Brownian setting) which can be proved by using potential theory.   

\subsection{Proof of theorem \ref{th:sin}}
First we state and prove a few results that we will use in the proof
of theorem \ref{th:sin}.
We remind the classical local limit theorem for the simple random walk
(cf. \cite{cf:La91}):
\begin{theorem}[Local limit theorem]\label{th:llt}
For $n \in \N$ and $x \in \Z^d$, we say that $n$ and $x$ have the same parity
and write $n \leftrightarrow x$ if $n+\sum_{k=1}^{d}x_{k}$ is
even and we define $\bar{q}^{(n)}(x)$ to be the gaussian approximation
of $q^{(n)}(x)$:         
$$\bar{q}^{(n)}(x)\overset{def.}{=}2(\frac{d}{2\pi
  n})^{\frac{d}{2}}\exp^{-\frac{d\mid x \mid^2}{2n}}.$$
With these notations, we have:

\begin{equation}
\label{eq:llt1}
  \sup_{ n \leftrightarrow x}\mid q^{(n)}(x)-\bar{q}^{(n)}(x)\mid = O(\frac{1}{n^{\frac{d}{2}+1}}). 
\end{equation}
In particular,
\begin{equation}
  \label{eq:llt2}
\sup_{ n \leftrightarrow
  x}q^{(n)}(x)=O(\frac{1}{n^{\frac{d}{2}}}).
\end{equation}
and if one fixes $A>0$, there exists $c>0$ such that
\begin{equation}
  \label{eq:llt3}
  \inftwo { n \leftrightarrow x}{\mid x \mid \leq A\sqrt{n}}q^{(n)}(x)
\geq c\frac{1}{n^{\frac{d}{2}}}.
\end{equation}

\end{theorem}

We will need the following obvious corollary of theorem \ref{th:llt}
wich can be understood as an absolute continuity result: 
\begin{corollary}\label{cor:moi}
Let $t \in ]0,1[$ and $A>0$. There exists a constant $C(A,d)>0$ such
that:
\begin{align} 
 \forall f \geq 0 \; \forall n \; \sup_{\mid y -x \mid \leq A\sqrt{n}}
&P^x \otimes P^x(f((\omega_{k},\tilde{\omega}_{k})_{k \leq \lfloor nt \rfloor} \mid
\omega_{n}=y,\tilde{\omega}_{n}=y) \nonumber \\ 
& \leq
\frac{C(A,d)}{(1-t)^{d}}P^x \otimes P^x(f((\omega_{k},\tilde{\omega}_{k})_{k \leq \lfloor
  nt \rfloor})). \label{eq:moi}
\end{align} 
\end{corollary}

\proof
By developping the left hand side of the inequality:

\begin{align*}
& P^x \otimes P^x(f((\omega_{k},\tilde{\omega}_{k})_{k \leq \lfloor nt
  \rfloor}) 
\mid \omega_{n}=y,\tilde{\omega}_{n}=y)\\ 
&=\sumtwo{z_{1}, \ldots, z_{\lfloor nt \rfloor} \in \Z^d} 
{\tilde{z}_{1}, \ldots, \tilde{z}_{\lfloor nt \rfloor} \in \Z^d}
q^{(1)}(z_{1}-x) \ldots q^{(1)}(z_{\lfloor nt \rfloor}-z_{\lfloor nt
  \rfloor-1}) \nonumber\\
& q^{(1)}(\tilde{z}_{1}-x) \ldots q^{(1)}(\tilde{z}_{\lfloor nt \rfloor}-\tilde{z}_{\lfloor nt
  \rfloor-1}) \\
&f(z_{1}, \ldots, z_{\lfloor nt \rfloor},\tilde{z}_{1}, \ldots, \tilde{z}_{\lfloor nt \rfloor})
\frac{ q^{(n-\lfloor nt \rfloor)}(y-z_{\lfloor nt
    \rfloor})}{q^{(n)}(y-x)} \\
& \frac{ q^{(n-\lfloor nt \rfloor)}(y-\tilde{z}_{\lfloor nt
    \rfloor})}{q^{(n)}(y-x)}. 
\end{align*}
By the local limit theorem \ref{th:llt},
\begin{equation*}
\frac{ q^{(n-\lfloor nt \rfloor)}(y-z_{\lfloor nt
    \rfloor})}{q^{(n)}(y-x)} \underset{(\ref{eq:llt2},\ref{eq:llt3})}{\leq} C' (\frac{n}{n-\lfloor nt
    \rfloor})^{\frac{d}{2}}\leq  \frac{C'}{(1-t)^{\frac{d}{2}}}.
\end{equation*}
Similarly,
\begin{equation*}
\frac{ q^{(n-\lfloor nt \rfloor)}(y-z_{\lfloor nt
    \rfloor})}{q^{(n)}(y-x)} \leq \frac{C'}{(1-t)^{\frac{d}{2}}}.
\end{equation*}
\qed

In order to prove theorem \ref{th:sin}, we will also need to use a result that
comes from discrete potential theory. For a complete overview of
potential theory for discrete Markov chains, we refer to \cite{cf:Wo}.

\begin{lemma}\label{lem:yo}
For $d \geq 3$ and $v: \Z^d \times \Z^d \longrightarrow \R$ a bounded function, define 

$$\Phi(x,y)=P^{x} \otimes P^{y}(e^{\sum_{k=1}^{\infty}v(\omega_{k},\tilde{\omega}_{k})}).$$
Suppose that
$$0<\inf_{x,y \in \Z^d}\Phi(x,y) \leq \sup_{x,y \in \Z^d}\Phi(x,y)< \infty.$$
Then there exists a constant $C \in ]0,\infty[$ such that 
\begin{equation}
\label{ine:yod} \sup_{x,y \in \Z^d}P^{x}\otimes P^{y}(e^{\sum_{k=1}^{n}v(\omega_{k},\tilde{\omega}_{k})}\mid
f(\omega_{n},\tilde{\omega}_{n})\mid)
\leq  \frac{C}{n^{d}}\sum_{x,y \in \Z^d}\mid f(x,y)\mid 
\end{equation}
for all $f$ in $L^{1}(\Z^{2d})$ and $n \geq 1$.
\end{lemma}

\proof
We will show inequality (\ref{ine:yod}) for $n$ even, the case $n$ odd
being similar.   
Let $(\omega_{n})_{n \geq 0}$ denote the simple random walk on
$\Z^d$. By theorem 4.18 in \cite{cf:Wo}, $(\omega_{2n})_{n \geq 0}$
satisfies the d-isoperimetric inequality ($IS_{d}$ therein) on its underlying graph. By
remark 4.11 in \cite{cf:Wo}, $(\omega_{2n},\tilde{\omega}_{2n})_{n
  \geq 0}$ satisfies the 2d-isoperimetric inequality on its underlying
graph. Consider the Markov chain in $\Z^d \times \Z^d$ with kernel:
\begin{equation*} 
K((x,y),(x',y'))=\sum_{(z,\tilde{z}) \in \Z^{2d}}\frac{1}{\Phi(x,y)}e^{v(z,\tilde{z})}e^{v(x',y')}p((x,y),(z,\tilde{z}))p((z,\tilde{z}),(x',y'))
\end{equation*}
where $p$ is the transition kernel of $(\omega_{n},\tilde{\omega}_{n})_{n
\geq 0}$. The transition kernel $K$ is reversible with invariant measure
$m(x,y)=(\Phi(x,y))^2e^{v(x,y)}$. By assumption, we have 
\begin{equation*}
0 < \inf_{x,y \in \Z^d}m(x,y) \leq
\sup_{x,y \in \Z^d}m(x,y) < \infty.  
\end{equation*}
By assumption, there exists $c,C>0$ such that for all $(x,y),(x',y')$ in $\Z^d$  
\begin{equation*}
 cp^{(2)}((x,y),(x',y'))\leq
K((x,y),(x',y')) \leq C p^{(2)}((x,y),(x',y'))
\end{equation*}
where $p^{(2)}$ is the transition kernel of $(\omega_{2n},\tilde{\omega}_{2n})_{n
\geq 0}$. Therefore $K$ satisfies the 2d-isoperimetric inequality on
its underlying graph. By corollary 14.5 in \cite{cf:Wo},  
\begin{equation*}
\sup_{x,y \in \Z^d}\frac{1}{\Phi(x,y)}P^{x}\otimes P^{y}(e^{\sum_{k=1}^{2n}v(\omega_{k},\tilde{\omega}_{k})}\mid
f(\omega_{2n},\tilde{\omega}_{2n})\mid \Phi(\omega_{2n},\tilde{\omega}_{2n}))
\leq C \frac{1}{n^{d}}\sum_{x,y \in \Z^d}\mid f(x,y)\mid
\end{equation*}
for all $f$ in $L^{1}(\Z^{2d})$ and $n \geq 1$. The inequality
(\ref{ine:yod}) follows by
using the boundedness of $v$ and the assumption on $\Phi$. \qed 

We can now state the following usefull corollary of lemma \ref{lem:yo}:
\begin{corollary}\label{cor:moi2}
Let $A>0$ and $x \in \Z^d$. Under the assumptions of lemma \ref{lem:yo}, there exists $C \in
]0,\infty[$ such that:
\begin{equation}
\label{eq:moi2}\forall n \; \; \sup_{\mid y -x \mid \leq
  A\sqrt{n}}P^{x} \otimes P^x(e^{\sum_{k=1}^{n}v(\omega_{k},\tilde{\omega}_{k})}\mid
\omega_{n}=y, \tilde{\omega}_{n}=y)
\leq C.
\end{equation}
\end{corollary}
\proof
Let $y \in \Z^d$ be such that $\mid y-x\mid\leq A\sqrt{n}$.  By
applying inequality (\ref{ine:yod}) with $f=1_{\lbrace y,y\rbrace}$, we get:
\begin{equation*}
P^{x}\otimes P^x(e^{\sum_{k=1}^{n}v(\omega_{k},\tilde{\omega}_{k})}\mid
\omega_{n}=y,\tilde{\omega}_{n}=y) \leq \frac{C}{n^{d}P^{x}\otimes P^x(\omega_{n}=y,\tilde{\omega}_{n}=y)}
\underset{(\ref{eq:llt3})}{\leq} C'.
\end{equation*} \qed

We can now prove theorem \ref{th:sin}.

\noindent\emph{Proof of theorem \ref{th:sin}.} Let $l_{n}$ be a sequence tending to infinity and such that $\forall n \;
l_{n} \leq n/2$.
First, we compare in $L^2$ the quantity $P^{x}(e_{1,n}\mid
\omega_{n}=y)$ with $P^{x}(e_{1,l_{n}}e_{n-l_{n},n} \mid
\omega_{n}=y)$. Therefore we compute:

\begin{align}
Q(P^{x}(e_{1,n}-e_{1,l_{n}}e_{n-l_{n},n} \mid
\omega_{n}=y))^2 &=P^x \otimes
P^x(e^{\lambda_{2}(\beta)N_{1,n}}-e^{\lambda_{2}(\beta)N_{1,l_{n}}}e^{\lambda_{2}(\beta)N_{n-l_{n},n}}\mid
\omega_{n}=y,\tilde{\omega}_{n}=y)  \nonumber \\
 &\leq P^x \otimes
P^x(e^{\lambda_{2}(\beta)N_{1,l_{n}}}e^{\lambda_{2}(\beta)N_{n-l_{n},n}}
\Delta_{n}\mid
\omega_{n}=y,\tilde{\omega}_{n}=y) \nonumber
\end{align}
 
with
$$\Delta_{n}=e^{\lambda_{2}(\beta)N_{l_{n},n-l_{n}}}-1.$$

Let $\delta>0$ be such that $(1+\delta)\lambda_{2}(\beta)<
\ln(\frac{1}{\pi_{d}})$. We remind that this implies:
\begin{align*}
\sup_{x,y \in \Z^d}P^{x} \otimes
P^{y}(e^{(1+\delta)\lambda_{2}(\beta)\sum_{k=1}^{\infty}1_{\omega_{k}=\tilde{\omega}_{k}}})
&=P^{x} \otimes
P^{x}(e^{(1+\delta)\lambda_{2}(\beta)\sum_{k=1}^{\infty}1_{\omega_{k}=\tilde{\omega}_{k}}})\\ 
&=P
\otimes P(e^{(1+\delta)\lambda_{2}(\beta)N_{1,\infty}}) < \infty. 
\end{align*}
Using inequality (\ref{eq:moi2}) with $v(x,y)=(1+\delta)\lambda_{2}(\beta)1_{x=y}$, there exists $C>0$ such that:
\begin{equation}\label{eq:co}
\sup_{n,\mid y-x\mid \leq A\sqrt{n}}P^x \otimes
P^x(e^{(1+\delta)\lambda_{2}(\beta)N_{1,n}} \mid
  \omega_{n}=y,\tilde{\omega}_{n}=y)\leq C.
\end{equation}  
Let $\epsilon, M$ be two
  positive numbers such that $e^{\lambda_{2}(\beta)\epsilon}-1< M$. By
  writing
  $$1=1_{\Delta_{n}<e^{\lambda_{2}(\beta)\epsilon}-1}+1_{M<\Delta_{n}}+1_{e^{\lambda_{2}(\beta)\epsilon}-1
    \leq \Delta_{n}\leq M},$$ we get

\begin{align*}
&P^x \otimes P^x(
e^{\lambda_{2}(\beta)(N_{1,l_{n}}+N_{n-l_{n},n})}\Delta_{n}\mid
  \omega_{n}=y,\tilde{\omega}_{n}=y)
\leq
C(e^{\lambda_{2}(\beta)\epsilon}-1)+\frac{C}{M^\delta} \nonumber \\ 
&+M P^x \otimes
P^x(1_{\Delta_{n}\geq
  e^{\lambda_{2}(\beta)\epsilon}-1}e^{\lambda_{2}(\beta)N_{1,n}}\mid
  \omega_{n}=y,\tilde{\omega}_{n}=y).
\end{align*}
Let $q>1$ be such that $\frac{1}{q}+\frac{1}{1+\delta}=1$. By Holder's
inequality and inequality (\ref{eq:moi2}), we get 
\begin{align*}
&P^x \otimes
P^x(1_{\Delta_{n}\geq
  e^{\lambda_{2}(\beta)\epsilon}-1}e^{\lambda_{2}(\beta)N_{1,n}}\mid
  \omega_{n}=y,\tilde{\omega}_{n}=y)) \\ 
&\leq (P^x \otimes
P^x(\Delta_{n}\geq
  e^{\lambda_{2}(\beta)\epsilon}-1\mid
  \omega_{n}=y,\tilde{\omega}_{n}=y))^{\frac{1}{q}}C^{\frac{1}{1+\delta}}.
\end{align*}
But, since $N_{l_{n},n-l_{n}}$ is integer valued, we get uniformly on $\mid y-x\mid \leq A\sqrt{n}$:  
\begin{align*}
&P^x \otimes
P^x(\Delta_{n}\geq
  e^{\lambda_{2}(\beta)\epsilon}-1\mid
  \omega_{n}=y,\tilde{\omega}_{n}=y) \\
&=P^x \otimes
P^x(N_{l_{n},n-l_{n}} \geq 1 \mid
  \omega_{n}=y,\tilde{\omega}_{n}=y) \\
&\leq P^x \otimes
P^x(N_{l_{n},n/2} \geq 1\mid
  \omega_{n}=y,\tilde{\omega}_{n}=y)) \\
&+P^x \otimes
P^x(N_{n/2,n-l_{n}} \geq 1\mid
  \omega_{n}=y,\tilde{\omega}_{n}=y)) \\
&= P^x \otimes
P^x(N_{l_{n},n/2} \geq 1\mid
  \omega_{n}=y,\tilde{\omega}_{n}=y)) \\
&+P^y \otimes
P^y(N_{l_{n},n/2} \geq 1\mid
  \omega_{n}=x,\tilde{\omega}_{n}=x)) \; \; (symmetry) \\
&\underset{(\ref{eq:moi})}{\leq} C'P^x \otimes
P^x(N_{l_{n},n/2} \geq 1) \\
&+C'P^y \otimes
P^y(N_{l_{n},n/2} \geq 1) \\
&= 2C'P\otimes P(N_{l_{n},n/2} \geq 1) \underset{n \to \infty}{\rightarrow}0.
\end{align*}
We have used in the limit above the fact that $N_{0,\infty}<\infty \; \; P\otimes
P-a.s.$ and that $l_{n}\underset{n \to \infty}{\longrightarrow} 0.$ Therefore, we get 
$$\varlimsup_{n \to \infty}\sup_{\mid y-x\mid \leq A\sqrt{n}}P^x \otimes P^x(
e^{\lambda_{2}(\beta)(N_{1,l_{n}}+N_{n-l_{n},n})}\Delta_{n}\mid
  \omega_{n}=y,\tilde{\omega}_{n}=y)
\leq
C(e^{\lambda_{2}(\beta)\epsilon}-1)+\frac{C}{M^\delta}.$$
We conclude that the above limit is equal to $0$ by letting $\epsilon \downarrow 0$ and $M \uparrow \infty$.

From now on, we suppose that $l_{n}=o(n^a)$ for some $a <
\frac{1}{2}$. By the Markov property of the simple random walk, we
get:
 
$$P^{x}(e_{1,l_{n}}e_{n-l_{n},n} \mid
\omega_{n}=y)=\sum_{\mid z_{1}-x \mid \leq l_{n},\mid y- z_{2} \mid \leq l_{n}} 
P^{x}(e_{1,l_{n}}1_{\omega_{l_{n}}=z_{1}})
\frac{q^{(n-2l_{n})}(z_{2}-z_{1})}{q^{(n)}(y-x)}
P^{z_{2}}(e_{n-l_{n},n}1_{\omega_{n}=y}).$$
By symmetry of the simple random walk, we have:
\begin{align*}
\sum_{\mid y- z_{2} \mid \leq
  l_{n}}P^{z_{2}}(e_{n-l_{n},n}1_{\omega_{n}=y}) &=
\sum_{\mid y- z_{2} \mid \leq
  l_{n}}P^{y}(\overset{\leftarrow}{e}_{1,l_{n}}1_{\omega_{n}=z_{2}}) \\
& =P^{y}(\overset{\leftarrow}{e}_{1,l_{n}}).
\end{align*}
Therefore,
\begin{align*}
&Q((P^{x}(e_{1,l_{n}}e_{n-l_{n},n} \mid
\omega_{n}=y)-P^{x}(e_{1,l_{n}})P^{y}(\overset{\leftarrow}{e}_{1,l_{n}}))^2)
 \qquad \qquad \\ 
&\qquad =\sumtwo{\mid z_{1}-x \mid \leq l_{n},\mid y- z_{2} \mid \leq
  l_{n}}{\mid z_{1}'-x \mid \leq l_{n},\mid y- z_{2}' \mid \leq
  l_{n}}\delta_{n}^{z_{1},z_{2},x,y}
\delta_{n}^{z_{1}',z_{2}',x,y} \times \\ 
&\qquad P^{x}\otimes P^{x}
(e^{\lambda^{2}(\beta)N_{1,l_{n}}}1_{\omega_{l_{n}}=z_{1}}
1_{\tilde{\omega}_{l_{n}}=z_{1}'}) P^{z_{2}}
\otimes
P^{z_{2}'}(e^{\lambda^{2}(\beta)N_{n-l_{n},n}}
1_{\omega_{n}=y}1_{\tilde{\omega}_{n}=y})
\end{align*}
where 
$$\delta_{n}^{z,w,x,y}=\frac{q^{(n-2l_{n})}(w-z)}{q^{(n)}(y-x)}-1.$$
The idea is that, by the classical local limit theorem, we get in the
previous sum the following estimate:
$$\frac{q^{(n-2l_{n})}(z_{2}-z_{1})}{q^{(n)}(y-x)}\approx
\frac{\bar{q}^{(n-2l_{n})}(z_{2}-z_{1})}{\bar{q}^{(n)}(y-x)} \approx 1.$$
Let us make this statement rigorous and obtain inequality
(\ref{ine:3}) below. We use the notations of
theorem \ref{th:llt}  and decompose
$\delta_{n}^{z,w,x,y}$ into three terms:
$$\delta_{n}^{z,w,x,y}=\delta_{1,n}^{z,w,x,y}+\delta_{2,n}^{z,w,x,y}+\delta_{3,n}^{z,w,x,y}$$
where 
$$\delta_{1,n}^{z,w,x,y}=\frac{q^{(n-2l_{n})}(w-z)-\bar{q}^{(n-2l_{n})}(w-z)}{q^{(n)}(y-x)},
\;
\delta_{2,n}^{z,w,x,y}=\frac{\bar{q}^{(n-2l_{n})}(w-z)-\bar{q}^{(n)}(y-x)}{q^{(n)}(y-x)},$$
$$\delta_{3,n}^{z,w,x,y}=\frac{\bar{q}^{(n)}(y-x)-q^{(n)}(y-x)}{q^{(n)}(y-x)}.$$ 
An application of (\ref{eq:llt1}) and (\ref{eq:llt3}) gives for $j=1,3$:  
$$\suptwo{\mid z-x \mid \leq l_{n},\mid y- w \mid \leq
  l_{n}}{\mid y -x \mid \leq A\sqrt{n}} \mid \delta_{j,n}^{z,w,x,y}\mid
= O(\frac{1}{n}) \underset{n \to
  \infty}{\longrightarrow} 0.$$

An application of (\ref{eq:llt3}) gives:
\begin{align*}
\mid \delta_{2,n}^{z,w,x,y}\mid &= \mid
\frac{\bar{q}^{(n-2l_{n})}(w-z)}{q^{(n)}(y-x)}\mid \mid 1-
\frac{\bar{q}^{(n)}(y-x)}{\bar{q}^{(n-2l_{n})}(w-z)} \mid \\
&\underset{(\ref{eq:llt3})}{\leq} C \mid 1-
\frac{\bar{q}^{(n)}(y-x)}{\bar{q}^{(n-2l_{n})}(w-z)} \mid. \\
&\leq C \mid 1- (\frac{n-2l_{n}}{n})^{\frac{d}{2}}e^{\frac{d\mid
    w-z\mid^2}{2(n-2l_{n})}-\frac{d\mid y-x\mid^2}{2n}} \mid. \\
\end{align*}
It is not hard to show that:
$$\suptwo{\mid z-x \mid \leq l_{n},\mid y- w \mid \leq
  l_{n}}{\mid y -x \mid \leq A\sqrt{n}}\mid 1- (\frac{n-2l_{n}}{n})^{\frac{d}{2}}e^{\frac{d\mid
    w-z\mid^2}{2(n-2l_{n})}-\frac{d\mid y-x\mid^2}{2n}} \mid \underset{n \to
  \infty}{\longrightarrow} 0$$
so we have 
$$\suptwo{\mid z-x \mid \leq l_{n},\mid y- w \mid \leq
  l_{n}}{\mid y -x \mid \leq A\sqrt{n}}\mid \delta_{2,n}^{z,w,x,y}\mid \underset{n \to
  \infty}{\longrightarrow} 0.$$

Finally, we get:
\begin{align}
\sup_{\mid y-x \mid \leq A\sqrt{n}}&Q((P^{x}(e_{1,l_{n}}e_{n-l_{n},n} \mid
\omega_{n}=y)-P^{x}(e_{1,l_{n}})P^{y}(\overset{\leftarrow}{e}_{1,l_{n}}))^2)
\nonumber\\
 &\leq \suptwo{\mid z-x \mid \leq l_{n},\mid y- w \mid \leq
  l_{n}}{\mid y -x \mid \leq A\sqrt{n}}\mid
\delta_{n}^{z,w,x,y}\mid^{2} (P\otimes P(e^{\lambda_{2}(\beta)(1+N_{1,l_{n}})}))^{2}\underset{n \to
  \infty}{\longrightarrow} 0. \label{ine:3}
\end{align}
Therefore, we get the expansion (\ref{eq:sin}). To get the expansion
(\ref{eq:sin2}), observe that
\begin{equation*}
P^{x}(e_{1,l_{n}})\overset{L^{2}(Q)}{\underset{n \to
    \infty}{\longrightarrow}}
Z_{\infty}^{x}
\end{equation*}
and, by symmetry, 
\begin{align*}
\sup_{y \in \Z^{d}}&Q((P^{y}(\overset{\leftarrow}{e}_{1,l_{n}})-
P^{y}(\overset{\leftarrow}{e}_{1,n}))^{2}) \\
&= P \otimes P(e^{\lambda_{2}(\beta)(1+N_{1,l_{n}-1})}-e^{\lambda_{2}(\beta)(1+N_{1,n-1})})\underset{n \to
    \infty}{\longrightarrow} 0.
\end{align*}

\qed

\subsection{Proof of theorem \ref{th:sin2}}

In order to prove theorem \ref{th:sin2}, we adapt in detail the
previous proof to the Brownian setting.
In the discrete setting, there are three key intermediate results: the
local limit theorem \ref{th:llt}, corollary \ref{cor:moi} and corollary
\ref{cor:moi2}.       
In the continuous setting, we do not need any local limit theorem since
Brownian motion is already a gaussian process. Therefore, we only
require a Brownian analogue to corollary \ref{cor:moi} and corollary
\ref{cor:moi2}.  
The following construction of the Brownian bridge can be found in the
appendix of \cite{cf:Szn98}:

\begin{proposition}
For $x,y \in \R^d$, $t>0$, there exists a unique probability measure
$P_{t}^{x,y}$ on $C([0,1],\R^d)$ such that for $s \in [0,t[, A \in \mathcal{F}_{s}$:
\begin{equation}
 \label{eq:szn}
 P_{t}^{x,y}(A)=\frac{1}{p(t,x,y)}P^{x}(1_{A}p(t-s,\omega_{s},y)) 
\end{equation}
$y \rightarrow P_{t}^{x,y}$ is a regular conditional probability of
$P_{t}^{x}$ given $\omega_{t}=y$.
\end{proposition}
In the sequel, we will always work with the representation
(\ref{eq:szn}) of Brownian bridge.
With this representation, we can now easily prove the brownian
analogue of corollary \ref{cor:moi}:

\begin{corollary}
Let $s \in ]0,1[$ and $A>0$. There exists a constant $C(A,d)>0$ such that 
\begin{align}
\forall f \geq 0 \; \; \forall t>0 \sup_{\mid y -x \mid \leq
  A\sqrt{t}}&P^{x}(f((\omega_{u})_{u \leq st})\mid \omega_{t}=y)
  \nonumber \\ 
&\leq 
\frac{C(A,d)}{(1-s)^{\frac{d}{2}}}P^{x}(f((\omega_{u})_{u \leq
  st})). \label{eq:moi'}
\end{align}
\end{corollary}
\proof
If $\mid y -x \mid \leq
  A\sqrt{t}$ then 

\begin{align*}
P^{x}(f((\omega_{u})_{u \leq st})\mid \omega_{t}=y) &=\frac{(2\pi
t)^{\frac{d}{2}}}{(2\pi
t(1-s))^{\frac{d}{2}}}e^{\frac{\mid y-x
\mid^2}{2t}}P^{x}(e^{-\frac{\mid y-\omega_{st}
\mid^2}{2t(1-s)}}f((\omega_{u})_{u \leq st})) \\
&\leq \frac{e^{A^2/2}}{(1-s)^{\frac{d}{2}}}P^{x}
(f((\omega_{u})_{u \leq st})).
\end{align*}
\qed

The Brownian analogue to lemma \ref{lem:yo} is a slight variation of
Lemma 3.1.3. in \cite{cf:CY04}.

\begin{lemma}
For $d \geq 3$ and $v:\R^d \longrightarrow \R$ a bounded, compactly
supported measurable function, define
$$ \Phi(x,y)=P^{x} \otimes P^{y}(e^{\int_{0}^{\infty}v(\tilde{\omega}_{s}-\omega_{s})ds}).$$
Suppose that 
$$0 < \inf_{x,y \in \R^{2d}}\Phi(x,y) \leq \sup_{x,y \in \R^{2d}}\Phi(x,y) <
\infty.$$
Then there exists a constant $C>0$ such that 
\begin{equation}
\label{ine:yo1} \sup_{x,y \in
  \R^{2d}}P^{x} \otimes P^{y}(e^{\int_{0}^{t}v(\tilde{\omega}_{s}-
\omega_{s})ds})\mid
f(\omega_{t},\tilde{\omega}_{t})\mid) \leq
\frac{C}{t^{d}}\int_{\R^{2d}}\mid f(x,y) \mid dxdy 
\end{equation}
for all $f$ in $L^1(\R^{2d})$ and $t>0$.
\end{lemma}

\proof
By using the same arguments than the ones in the proof of Lemma
3.1.3. in \cite{cf:CY04}, all we have to prove is 
\begin{equation*} 
\forall F \in C_{c}^{\infty}(\R^{2d}) \; \; \int_{\R^{2d}}(\frac{1}{2}\bigtriangledown_{x,y}F.\bigtriangledown_{x,y}\Phi-v(y-x)F(x,y)\Phi(x,y))dxdy=0.
\end{equation*}
Since $(\tilde{\omega}_{s/2}-\omega_{s/2})_{s \geq 0}$ is a brownian
motion, we have that $\Phi(x,y)=\tilde{\Phi}(y-x)$ where:
\begin{equation*}
\forall z \in \R^d \; \;
\tilde{\Phi}(z)=P^z(e^{\int_{0}^{\infty}\frac{1}{2}v(\omega_{s})ds}).
\end{equation*}
By equation (3.19) in the proof Lemma 3.1.3. in \cite{cf:CY04}, we
have:
\begin{equation} \label{eq:yoyo}
\forall g \in C_{c}^{\infty}(\R^{d}) \; \; \int_{\R^{d}}(\bigtriangledown_{\tilde{y}}g(\tilde{y}).\bigtriangledown_{\tilde{y}}\tilde{\Phi}-v(\tilde{y})g(\tilde{y})\tilde{\Phi}(\tilde{y}))d\tilde{y}=0.
\end{equation}
By making the change of variable $(\tilde{x},\tilde{y})=(y+x,y-x)$, $F(x,y)=f(\tilde{x},\tilde{y})$ we get:
\begin{equation*}
\bigtriangledown_{x,y}F.\bigtriangledown_{x,y}\Phi=-(\bigtriangledown_{\tilde{x}}f-
\bigtriangledown_{\tilde{y}}f)\bigtriangledown_{\tilde{y}}\tilde{\Phi}
+(\bigtriangledown_{\tilde{x}}f+\bigtriangledown_{\tilde{y}}f)
\bigtriangledown_{\tilde{y}}\tilde{\Phi}
=2\bigtriangledown_{\tilde{y}}f\bigtriangledown_{\tilde{y}}\tilde{\Phi}
\end{equation*}
Therefore, 
\begin{align*}
\forall F \in C_{c}^{\infty}(\R^{2d}) \; \;
&\int_{\R^{2d}}(\frac{1}{2}\bigtriangledown_{x,y}F.\bigtriangledown_{x,y}\Phi-v(y-x)f(x,y)\Phi(x,y))dxdy
\\&=\frac{1}{2^d}\int_{\R^{2d}}(\bigtriangledown_{\tilde{y}}f(\tilde{x},\tilde{y})\bigtriangledown_{\tilde{y}}\tilde{\Phi}(\tilde{y})-v(\tilde{y})f(\tilde{x},\tilde{y})
\tilde{\Phi}(\tilde{y}))d\tilde{x}d\tilde{y} \\
&=\frac{1}{2^d}\int_{\R^{d}}(\int_{\R^{d}}(\bigtriangledown_{\tilde{y}}f(\tilde{x},\tilde{y})\bigtriangledown_{\tilde{y}}\tilde{\Phi}(\tilde{y})-v(\tilde{y})f(\tilde{x},\tilde{y})
\tilde{\Phi}(\tilde{y}))d\tilde{y})d\tilde{x} \\
&\underset{(\ref{eq:yoyo})}{=}0
\end{align*}
\qed

We can now state the following analogue to corollary \ref{cor:moi2}: 

\begin{corollary}
Let $A >0$ and $x,y \in \R^d$. Under the above assumptions, there exists
$C > 0$ such that:
\begin{equation}
\forall t \; \; \sup_{\mid y - x \mid \leq A
  \sqrt{t}}P^x \otimes
  P^{x}(e^{\int_{0}^{t}v(\tilde{\omega}_{s}-\omega_{s})ds} \mid
  \omega_{t}=y,\tilde{\omega}_{t}=y) \leq C. \label{eq:moi2'}
\end{equation}
\end{corollary}

\proof
Let $r>0$ and $y \in \R^d$ such that  $\mid y - x \mid \leq
A\sqrt{t}$. By applying (\ref{ine:yo1}) with $f=1_{B((y,y),r)}$, we get:
$$P^x \otimes P^{x}(e^{\int_{0}^{t}v(\tilde{\omega}_{s}-\omega_{s})ds}1_{B((y,y),r)}(\omega_{t},\tilde{\omega}_{t}))\leq
\frac{C}{t^{d}}\mid B((y,y),r)\mid.$$
Therefore, 
\begin{align}
P^x\otimes
P^{x}(e^{\int_{0}^{t}v(\tilde{\omega}_{s}-\omega_{s})ds}1_{B((y,y),r)}(\omega_{t},\tilde{\omega}_{t}))&/
P^x\otimes P^{x}((\omega_{t},\tilde{\omega}_{t}) \in B((y,y),r))
\nonumber \\ 
&\leq
\frac{C}{t^{d}}\frac{\mid B((y,y),r)\mid}{P^x\otimes
  P^{x}((\omega_{t},\tilde{\omega}_{t}) \in
  B((y,y),r))}. \label{eq:truc}
\end{align}
As $r\downarrow 0$, a classical result on brownian bridges asserts
that the left handside of (\ref{eq:truc}) tends to  
$$P^x\otimes P^{x}(e^{\int_{0}^{t}v(\tilde{\omega}_{s}-\omega_{s})ds} \mid \omega_{t}=y,\tilde{\omega}_{t}=y).$$
As $r\downarrow 0$, the right handside of (\ref{eq:truc}) tends to
$$C\frac{e^{\frac{\mid y-x\mid^2}{t}}(2\pi
  t)^{d}}{t^{d}}\leq (2\pi)^{d}Ce^{A^2}.$$
\qed

\noindent\emph{Proof of theorem \ref{th:sin2}.}
The proof of theorem \ref{th:sin2} is quite similar but even simpler
than the proof of theorem
\ref{th:sin} since brownian motion is already gaussian. We will not
repeat the details but we indicate the main steps for convenience. Suppose that $\beta$ is such
that
\begin{equation*} 
\lambda_{2}(\beta)<\lambda(d).
\end{equation*}
There exists $\delta>0$ such that
$(1+\delta)\lambda_{2}(\beta)<\lambda(d)$. Using inequality
(\ref{eq:moi2'}) applied to $v(y-x)=(1+\delta)\lambda_{2}(\beta)\mid
U(y-x) \cap U(0)\mid$, we get the following analogue to (\ref{eq:co}): there
exists $C>0$ such that
\begin{equation}\label{eq:co'}
\sup_{t, \mid y-x\mid \leq A\sqrt{t}}P^{x}\otimes
P^{x}(e^{(1+\delta)\lambda_{2}(\beta)N_{0,t}}\mid
\omega_{t}=y,\tilde{\omega}_{t}=y)\leq C.
\end{equation}
Using inequality (\ref{eq:co'}) and inequality (\ref{eq:moi'}), we get 
\begin{equation*}
P^{x}(e_{0,t}\mid \omega_{t}=y) \approx P^{x}(e_{0,l_{t}}e_{t-l_{t},t}\mid \omega_{t}=y).
\end{equation*}
Using the Markov property and the symmetry of Brownian motion, we get
\begin{equation*}
P^{x}(e_{0,l_{t}}e_{t-l_{t},t}\mid \omega_{t}=y) \approx P^{x}(e_{0,l_{t}})P^{y}(\overset{\leftarrow}{e}_{0,l_{t}}).
\end{equation*}
\qed

\textbf{Acknowledgements}: I would like to thank my Ph.D. supervisor
Francis Comets for his help and suggestions.


\end{document}